\documentclass[12pt]{amsart}
\usepackage{amssymb}
\usepackage{graphicx}
\newtheorem*{thm}{Theorem}
\title[Poincar\'e-Birkhoff-Witt]
{A geometric proof of the Poincar\'e-Birkhoff-Witt Theorem}
\author[Michael Eastwood]{Michael Eastwood}
\address{\hskip-\parindent
School of Mathematical Sciences\\
University of Adelaide\\ 
SA 5005\\ 
Australia}
\email{meastwoo@member.ams.org}

\thanks{This work was also supported by the Simons Foundation grant 346300 and
the Polish Government MNiSW 2015--2019 matching fund. It was written whilst the
author was at the Banach Centre at IMPAN in Warsaw for the Simons Semester
`Symmetry and Geometric Structures.' This proof was presented in the author's
lectures on `Invariant Differential Operators' at the start of this Simons 
Semester.}

\begin{document}\raggedbottom
\begin{abstract}\small 
We use that the $n$-sphere for $n\geq 2$ is simply-connected to prove
the Poincar\'e-Birkhoff-Witt Theorem.
\end{abstract}    
\maketitle 

There are several equivalent statements of the
Poincar\'e-Birkhoff-Witt Theorem. The version we shall prove is as follows. 
\begin{thm}
Let ${\mathfrak{g}}$ be a Lie algebra. Define an equivalence relation on the 
tensor algebra $\bigotimes{\mathfrak{g}}$ by imposing the relations that
\begin{equation}\label{relations}\tag{$\star\star\star$}
a\otimes b-b\otimes a=[a,b]\end{equation}
as a two-sided ideal in $\bigotimes{\mathfrak{g}}$. Write the resulting
associative algebra as ${\mathfrak{U}}({\mathfrak{g}})$ and write $ab\cdots d$
for the equivalence class of $a\otimes b\otimes\cdots\otimes d$. Pick a basis
for ${\mathfrak{g}}$ and declare that an element $ab\cdots
d\in{\mathfrak{U}}({\mathfrak{g}})$ is in `canonical form' if and only if
$a,b,\ldots,d$ are basis elements with $a\leq b\leq\cdots\leq d$ with respect
to the ordering of the basis. Then elements in ${\mathfrak{U}}({\mathfrak{g}})$
may be consistently and uniquely written as
linear combinations of elements in canonical form. 
\end{thm}

An algebraic proof may be found, for example, in~\cite{J}. The rest of this 
article is devoted to a geometric proof.

To understand what the Poincar\'e-Birkhoff-Witt Theorem says, let us consider 
the case of three elements $a,b,c\in{\mathfrak{g}}$, which we suppose are 
basis elements in this order $a\leq b\leq c$, and that we would like to 
rewrite the element $cba\in{\mathfrak{U}}({\mathfrak{g}})$ (given in the 
`wrong' order) as a linear combination of canonically ordered elements. 
Certainly, we can use the equivalence relation (\ref{relations}) to try to 
reorder this element:
$$\begin{array}{rcl}cba&=&cab-c[a,b]\\
                          &=&acb-[a,c]b-c[a,b]\\
                          &=&abc-a[b,c]-[a,c]b-c[a,b],\end{array}$$
where we have firstly swopped $b$ and $a$ (and then followed our noses). The 
only problem is that one can firstly swop $c$ and $b$ instead:
$$\begin{array}{rcl}cba&=&bca-[b,c]a\\
                       &=&bac-b[a,c]-[b,c]a\\
                       &=&abc-[a,b]c-b[a,c]-[b,c]a,\end{array}$$
which is consistent if and only if the `second order' remainder terms agree:
$$a[b,c]+[a,c]b+c[a,b]=[a,b]c+b[a,c]+[b,c]a.$$
Fortunately, this is exactly the Jacobi identity:
$$[a,[b,c]]+[[a,c],b]+[c,[a,b]]=0.$$
We may arrange these calculations on a circle:
$$\begin{picture}(120,120)
\put(-80,60){\makebox(0,0){Figure~1}}
\put(60,60){\makebox(0,0){\large$\bullet$}}
\qbezier(60,60)(60,60)(60,100)
\qbezier(60,110)(60,110)(60,120)
\qbezier(60,60)(60,60)(94.64101616,80)
\qbezier(103.3012702,85)(103.3012702,85)(111.9615242,90)
\qbezier(60,60)(60,60)(94.64101616,40)
\qbezier(103.3012702,35)(103.3012702,35)(111.9615242,30)
\qbezier(60,60)(60,60)(25.35898384,80)
\qbezier(16.69872980,35)(16.69872980,35)(8.03847576,30)
\qbezier(60,60)(60,60)(25.35898384,40)
\qbezier(16.69872980,85)(16.69872980,85)(8.03847576,90)
\qbezier(60,60)(60,60)(60,20)
\qbezier(60,10)(60,10)(60,0)
\qbezier(50,104)(50,104)(70,104)
\qbezier(50,107)(50,107)(70,107)
\qbezier(103.1051178,73.33974596)(103.1051178,73.33974596)
        (93.10511778,90.66025404)
\qbezier(105.7031940,74.83974596)(105.7031940,74.83974596)
        (95.70319399,92.16025404)
\qbezier(93.10511778,29.33974596)(93.10511778,29.33974596)
        (103.1051178,46.66025404)
\qbezier(95.70319399,27.83974596)(95.70319399,27.83974596)
        (105.7031940,45.16025404)
\qbezier(50,16)(50,16)(70,16)
\qbezier(50,13)(50,13)(70,13)
\qbezier(16.89488222,46.66025404)(16.89488222,46.66025404)
        (26.89488222,29.33974596)
\qbezier(14.29680601,45.16025404)(14.29680601,45.16025404)
        (24.29680601,27.83974596)
\qbezier(26.89488222,90.66025404)(26.89488222,90.66025404)
        (16.89488222,73.33974596)
\qbezier(24.29680601,92.16025404)(24.29680601,92.16025404)
        (14.29680601,74.83974596)
\put(8,60){\makebox(0,0){$cba$}}
\put(32,99){\makebox(0,0){$bca+\cdots$}}
\put(95,99){\makebox(0,0){$bac+\cdots$}}
\put(122,60){\makebox(0,0){$abc+\cdots$}}
\put(32,22){\makebox(0,0){$cab+\cdots$}}
\put(95,22){\makebox(0,0){$acb+\cdots$}}
\end{picture}$$
where $\,\cdots$ denotes second order terms. Otherwise said, the Jacobi
identity is exactly what is needed so that an excursion through the symmetric
group ${\mathfrak{S}}_3$ on three letters
$$abc\rightsquigarrow bac\rightsquigarrow bca\rightsquigarrow
cba\rightsquigarrow cab\rightsquigarrow acb\rightsquigarrow abc$$
is consistent in~${\mathfrak{U}}({\mathfrak{g}})$. One can think of this as 
saying that there is no `holonomy' around the circle depicted in Figure~1. 

If we attempt a similar proof for four basis element $a\leq b\leq c\leq d$, 
then we run into trouble because there is no `follow your nose' method for 
reordering elements of the symmetric group~${\mathfrak{S}}_4$. Instead, we may 
picture ${\mathfrak{S}}_4$ as 24 countries in the plane arranged like this:
$$\begin{picture}(0,350)
\put(-155,260){\makebox(0,0){Figure~2}}
\put(-200,-110){{\includegraphics[width=14cm]{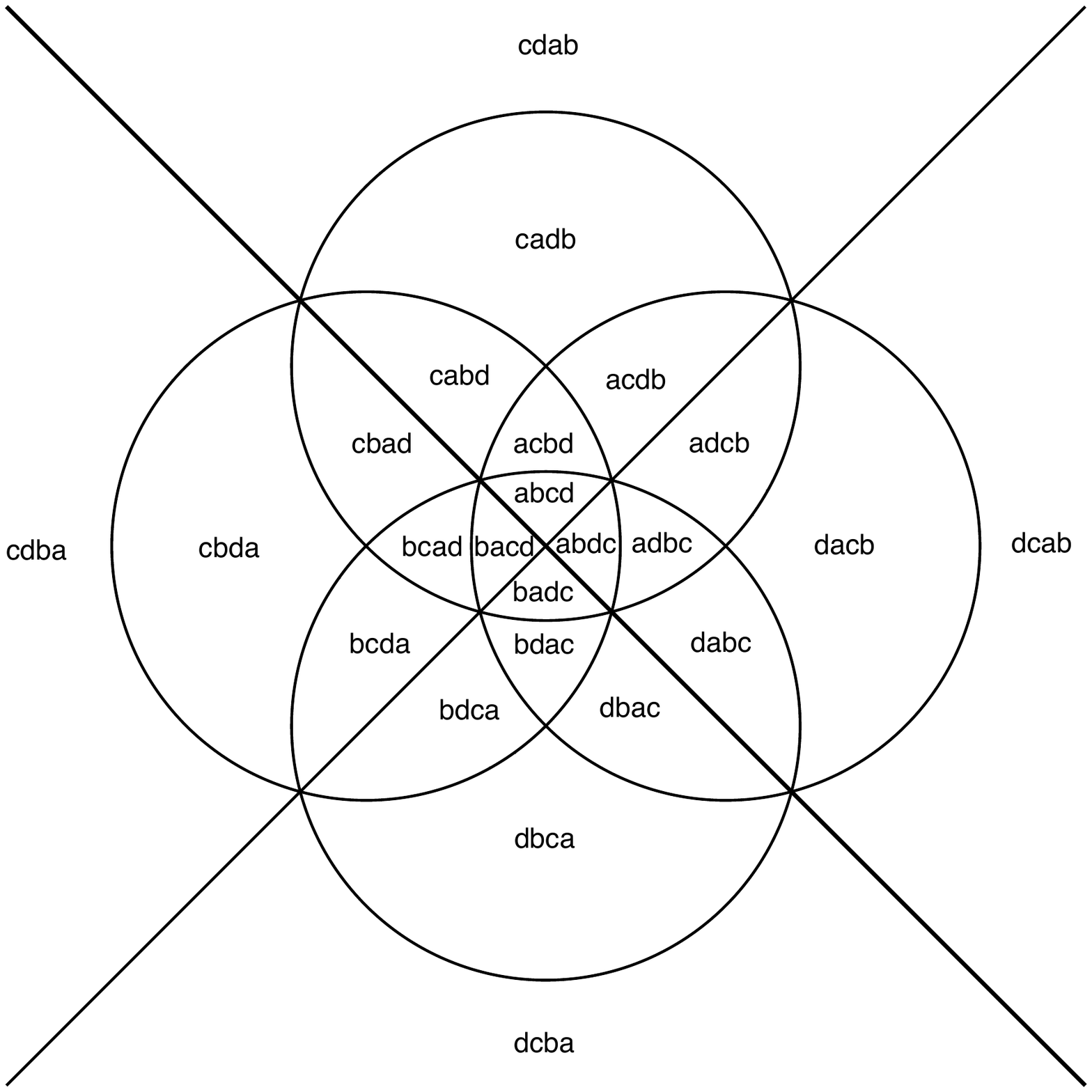}}}
\put(21,193){\makebox(0,0){\Large$\bullet$}}
\put(-21,193){\makebox(0,0){\Large$\bullet$}}
\put(21,151){\makebox(0,0){\Large$\bullet$}}
\put(-21,151){\makebox(0,0){\Large$\bullet$}}
\put(78.5,250.5){\makebox(0,0){\Large$\bullet$}}
\put(-78.5,250.5){\makebox(0,0){\Large$\bullet$}}
\put(78.5,93.5){\makebox(0,0){\Large$\bullet$}}
\put(-78.5,93.5){\makebox(0,0){\Large$\bullet$}}
\put(0,180){\makebox(0,0){$\bullet$}}
\thicklines
\qbezier(0,180)(10,180)(30,190)
\qbezier(30,190)(40,195)(50,240)
\qbezier(50,240)(60,285)(80,285)
\qbezier(80,285)(100,285)(110,200)
\qbezier(110,200)(120,115)(100,80)
\qbezier(100,80)(80,45)(-120,65)
\qbezier(-120,65)(-140,67)(-140,77)
\qbezier(-140,77)(-140,90)(-60,180)
\qbezier(-60,180)(-52,189)(0,180)
\put(53,260){\rotatebox{68}{$\blacktriangleright$}}
\put(109.4,140){\rotatebox{-90}{$\blacktriangleright$}}
\put(-87.4,150){\rotatebox{45}{$\blacktriangleright$}}
\end{picture}$$
Also depicted is a typical excursion through ${\mathfrak{S}}_4$
starting and finishing at $abcd$, namely
$$\addtolength{\arraycolsep}{-3pt}\begin{array}{ccccccccccccccc}abcd
&\rightsquigarrow& abdc
&\rightsquigarrow& adbc
&\rightsquigarrow& adcb
&\rightsquigarrow& acdb
&\rightsquigarrow& cadb
&\rightsquigarrow& cdab
&\rightsquigarrow& dcab\\
\rotatebox[origin=c]{90}{$\rightsquigarrow$}&&&&&&&&&&&&&&
\rotatebox[origin=c]{270}{$\rightsquigarrow$}\\
bacd&&&&&&&&&&&&&&dacb\\
\rotatebox[origin=c]{90}{$\rightsquigarrow$}&&&&&&&&&&&&&&
\rotatebox[origin=c]{270}{$\rightsquigarrow$}\\
bcad
&\raisebox{5pt}{\,\rotatebox{180}{$\rightsquigarrow$}\,}&cbad
&\raisebox{5pt}{\,\rotatebox{180}{$\rightsquigarrow$}\,}&cbda
&\raisebox{5pt}{\,\rotatebox{180}{$\rightsquigarrow$}\,}&cdba
&\raisebox{5pt}{\,\rotatebox{180}{$\rightsquigarrow$}\,}&dcba
&\raisebox{5pt}{\,\rotatebox{180}{$\rightsquigarrow$}\,}&dbca
&\raisebox{5pt}{\,\rotatebox{180}{$\rightsquigarrow$}\,}&dcba
&\raisebox{5pt}{\,\rotatebox{180}{$\rightsquigarrow$}\,}&dcab
\end{array}$$
We would like to see that this excursion is consistent. There are just $8$
points in the plane where $6$ countries come together. For example:
$$\begin{picture}(100,100)
\put(-100,40){\makebox(0,0){Figure 3}}
\put(50,50){\makebox(0,0){\large{$\bullet$}}}
\thicklines 
\qbezier(50,50)(50,50)(99.00332889,40.06653346)
\qbezier(50,50)(50,50)(83.10429882,87.47139442)
\qbezier(50,50)(50,50)(34.10096992,97.40486096)
\qbezier(50,50)(50,50)(0.99667111,59.93346654)
\qbezier(50,50)(50,50)(16.89570118,12.52860558)
\qbezier(50,50)(50,50)(65.89903008,2.59513904)
\put(40,0){\makebox(0,0){$adcb$}}
\put(0,40){\makebox(0,0){$acdb$}}
\put(0,80){\makebox(0,0){$cadb$}}
\put(50,95){\makebox(0,0){$cdab$}}
\put(100,70){\makebox(0,0){$dcab$}}
\put(100,20){\makebox(0,0){$dacb$}}
\qbezier(20,0)(20,0)(25,40)
\qbezier(25,40)(30,80)(60,80)
\qbezier(60,80)(80,80)(85,40)
\qbezier(85,40)(87,10)(85,0)
\put(29,62){\rotatebox{55}{$\blacktriangleright$}}
\end{picture}$$
These are the eight points where countries of the form 
$$a{*}{*}{*}\enskip\mbox{or}\enskip 
b{*}{*}{*}\enskip\mbox{or}\enskip
c{*}{*}{*}\enskip\mbox{or}\enskip
d{*}{*}{*}\enskip\mbox{or}\enskip
{*}{*}{*}a\enskip\mbox{or}\enskip
{*}{*}{*}b\enskip\mbox{or}\enskip
{*}{*}{*}c\enskip\mbox{or}\enskip
{*}{*}{*}d$$
meet at a vertex and these are marked by {\large$\bullet$} in Figure~2. 
The picture above is of the vertex ${*}{*}{*}b$ and one
recognises the circle from Figure~1 save that the elements $a,b,c$ have been
relabelled $d,c,a$. We saw earlier that this circle corresponds to a consistent
identity for three elements in~${\mathfrak{U}}({\mathfrak{g}})$ and now we
obtain a consistent identity for four elements in which $b$ simply goes along
for the ride. Geometrically, it means we may replace the path in Figure~3 by
$$\begin{picture}(100,100)
\put(50,50){\makebox(0,0){\large{$\bullet$}}}
\thicklines 
\qbezier(50,50)(50,50)(99.00332889,40.06653346)
\qbezier(50,50)(50,50)(83.10429882,87.47139442)
\qbezier(50,50)(50,50)(34.10096992,97.40486096)
\qbezier(50,50)(50,50)(0.99667111,59.93346654)
\qbezier(50,50)(50,50)(16.89570118,12.52860558)
\qbezier(50,50)(50,50)(65.89903008,2.59513904)
\put(40,0){\makebox(0,0){$adcb$}}
\put(0,40){\makebox(0,0){$acdb$}}
\put(0,80){\makebox(0,0){$cadb$}}
\put(50,95){\makebox(0,0){$cdab$}}
\put(100,70){\makebox(0,0){$dcab$}}
\put(100,20){\makebox(0,0){$dacb$}}
\qbezier(20,0)(22,10)(40,30)
\qbezier(40,30)(58,50)(70,35)
\qbezier(70,35)(82,20)(85,0)
\put(31.4,21){\rotatebox{50}{$\blacktriangleright$}}
\end{picture}$$
to obtain an alternative but simpler excursion through ${\mathfrak{S}}_4$,
which is consistent if and only if the original excursion is consistent. If we
can similarly pull paths through the other $5$ vertices where just four
countries come together, then we can reduce any excursion through
${\mathfrak{S}}_4$ to the trivial excursion (by a series of `simple jerks' in
the terminology of~\cite{G}) and our proof is complete. A typical example is
$$\begin{picture}(80,80)\thicklines
\put(-100,40){\makebox(0,0){Figure~4}}
\qbezier(0,0)(0,0)(80,80)
\qbezier(80,0)(80,0)(0,80)
\put(40,5){\makebox(0,0){$bcda$}}
\put(5,40){\makebox(0,0){$cbda$}}
\put(40,75){\makebox(0,0){$cbad$}}
\put(75,40){\makebox(0,0){$bcad$}}
\qbezier(0,10)(0,10)(30,40)
\qbezier(30,40)(60,70)(90,65)
\put(7.4,18){\rotatebox{45}{$\blacktriangleright$}}
\end{picture}$$
but vertices like this evidently have consistent holonomy
$$\addtolength{\arraycolsep}{-5pt}\begin{array}{rcccl}
&&cbad-cb[a,d]&&\\[-2pt]
&\rotatebox[origin=c]{45}{$=$}&&\rotatebox[origin=c]{-45}{$=$}\\[-2pt]
cbda&&&&bcad-[b,c]ad-bc[a,d]+[b,c][a,d]\\[-2pt]
&\rotatebox[origin=c]{-45}{$=$}&&\rotatebox[origin=c]{45}{$=$}\\[-2pt]
&&bcda-[b,c]da
\end{array}$$
without using the Jacobi identity. It is because we are transposing the first
two and the last two of four letters, and such transpositions commute
in~${\mathfrak{S}}_4$.

So now, we may consistently reorder any four elements in
${\mathfrak{U}}({\mathfrak{g}})$ and we ask about five elements and so on. 
We need a similar picture of the symmetric groups~${\mathfrak{S}}_N$
for all $N\geq 4$. To obtain such a picture, we now admit that Figure~2 was
obtained from a tessellation of the $2$-sphere by $24$ geodesic triangles with
angles $(\pi/2,\pi/3,\pi/3)$. Specifically, it was obtained by stereographic
projection so that great circles on the sphere are mapped to circles or
straight lines on the plane whilst angles are preserved.
$$\begin{picture}(200,200)
\put(0,0){{\includegraphics[width=7cm]{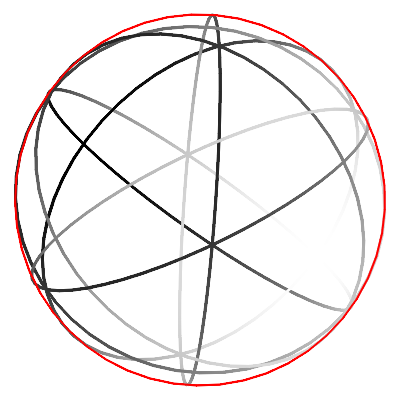}}}
\end{picture}$$
Therefore, a better
viewpoint on Figure~2 is as a triangulation of the $2$-sphere. {From} this
point of view there is one more `easy vertex,' as in Figure~4, out at infinity.
The fact that one can contract any excursion in ${\mathfrak{S}}_4$ to the
trivial excursion is due to there being no obstructions
\begin{itemize}
\item at the $6$ `easy vertices' (commuting transpositions),
\item at the $8$ `tricky vertices' (from the ${\mathfrak{S}}_3$ case),
\end{itemize} 
and the fact that the $2$-sphere is simply-connected. This triangulation of the
$2$-sphere is well-known in a different guise. It is obtained by letting the
Weyl group of the $A_3$ root system act on~${\mathbb{R}}^3$, as described, for
example, in~\cite{H}. The triangulation is obtained by intersecting the $24$
Weyl chambers with the unit sphere in~${\mathbb{R}}^3$. Since the Weyl group of
$A_3$ may be identified with ${\mathfrak{S}}_4$, one can pick a triangle to be
called the `fundamental triangle' and use the Weyl group action to identify any
element of ${\mathfrak{S}}_4$ with the triangle obtained as the corresponding
image of the fundamental triangle. This is how Figure~2 was obtained.

It is evident how to extend this to ${\mathfrak{S}}_N$ for all $N\geq 4$ and,
for the general pattern, it suffices to make sure that ${\mathfrak{S}}_5$
behaves as it should. The corresponding tessellation of the unit $3$-sphere is
by $120$ tetrahedra with dihedral angles
$(\pi/2,\pi/2,\pi/2,\pi/3,\pi/3,\pi/3)$ (each dihedral angle corresponds to a
pair of vertices from the Dynkin diagram $\begin{picture}(55,10)
\put(5,3){\line(1,0){45}} \put(5,3){\makebox(0,0){$\bullet$}}
\put(20,3){\makebox(0,0){$\bullet$}} \put(35,3){\makebox(0,0){$\bullet$}}
\put(50,3){\makebox(0,0){$\bullet$}} \end{picture}$, which are either adjacent
(angle $\pi/3$) or not (angle $\pi/2$)). To use the simple connectivity of the
$3$-sphere it now suffices to be able to move a path on the $3$-sphere 
through any edge of this tessellation.

As on the $2$-sphere, there are two cases. Firstly, there are the `easy edges,'
where just $4$ tetrahedra meet at right angles. On the Dynkin diagram, edges of
this type correspond to striking out all but two $2$ non-adjacent nodes
$$\begin{picture}(55,10)
\put(5,3){\line(1,0){45}} \put(5,3){\makebox(0,0){$\times$}}
\put(20,3){\makebox(0,0){$\bullet$}} \put(35,3){\makebox(0,0){$\times$}}
\put(50,3){\makebox(0,0){$\bullet$}} \end{picture}\qquad\begin{picture}(55,10)
\put(5,3){\line(1,0){45}} \put(5,3){\makebox(0,0){$\bullet$}}
\put(20,3){\makebox(0,0){$\times$}} \put(35,3){\makebox(0,0){$\times$}}
\put(50,3){\makebox(0,0){$\bullet$}} \end{picture}\qquad\begin{picture}(55,10)
\put(5,3){\line(1,0){45}} \put(5,3){\makebox(0,0){$\bullet$}}
\put(20,3){\makebox(0,0){$\times$}} \put(35,3){\makebox(0,0){$\bullet$}}
\put(50,3){\makebox(0,0){$\times$}} \end{picture},$$
in effect leaving the Weyl group of $A_1\times A_1$ as in Figure~$4$. It is
just the Abelian group ${\mathbb{Z}}_2\times{\mathbb{Z}}_2$. The `tricky edges'
are when $6$ tetrahedron meet at angle $\pi/3$. Tricky edges may be recorded on
the Dynkin diagram by striking out all but two adjacent nodes
$$\begin{tabular}{rcl}
\begin{picture}(55,10)
\put(5,3){\line(1,0){45}} \put(5,3){\makebox(0,0){$\times$}}
\put(20,3){\makebox(0,0){$\times$}} \put(35,3){\makebox(0,0){$\bullet$}}
\put(50,3){\makebox(0,0){$\bullet$}} \end{picture}&$\longleftrightarrow$
&permuting $ab{*}{*}{*}$ with $a,b$ held fixed,\\[5pt]
\begin{picture}(55,10)
\put(5,3){\line(1,0){45}} \put(5,3){\makebox(0,0){$\times$}}
\put(20,3){\makebox(0,0){$\bullet$}} \put(35,3){\makebox(0,0){$\bullet$}}
\put(50,3){\makebox(0,0){$\times$}} \end{picture}&$\longleftrightarrow$
&permuting $a{*}{*}{*}b$ with $a,b$ held fixed,\\[5pt]
\begin{picture}(55,10)
\put(5,3){\line(1,0){45}} \put(5,3){\makebox(0,0){$\bullet$}}
\put(20,3){\makebox(0,0){$\bullet$}} \put(35,3){\makebox(0,0){$\times$}}
\put(50,3){\makebox(0,0){$\times$}}
\end{picture}&$\longleftrightarrow$
&permuting ${*}{*}{*}ab$ with $a,b$ held fixed.
\end{tabular}$$
The tricky edges are not obstructed since the previous reasoning using the
Jacobi identity applies (notice that we are left with
$A_2=\begin{picture}(25,10) \put(5,3){\line(1,0){15}}
\put(5,3){\makebox(0,0){$\bullet$}} \put(20,3){\makebox(0,0){$\bullet$}}
\end{picture}$ and the Weyl group of $A_2$ is~${\mathfrak{S}}_3$). Looking
back, we see that Figure~1 is the root diagram for $A_2$.

\end{document}